\documentclass[12pt]{amsart}
\usepackage{enumitem}

\usepackage{amssymb, mathrsfs,mathcomp,amsfonts,longtable, hhline, textcomp, upgreek, graphicx, url}
\usepackage[matrix,arrow,curve]{xy}

\usepackage[OT2,T1]{fontenc}
\DeclareSymbolFont{cyrletters}{OT2}{wncyr}{m}{n}
\DeclareMathSymbol{\Sha}{\mathalpha}{cyrletters}{"58}

\newcommand{\B}{{\mathbf B}}
\newcommand{\CC}{\mathbb{C}}

\newcommand{\QQ}{\mathbb{Q}}

\newcommand{\ZZ}{\mathbb{Z}}
\newcommand{\PP}{\mathbb{P}}
\newcommand{\OOO}{{\mathscr{O}}} 
 
\newcommand{\MMM}{{\mathscr{M}}} 
 
\newcommand{\LLL}{{\mathscr{L}}}
\newcommand{\mumu}{{\boldsymbol{\mu}}}

\newcommand{\qW}{\operatorname{q}_{\operatorname{W}}}
\newcommand{\qQ}{\operatorname{q_{\QQ}}}
\newcommand{\DP}{\operatorname{DP}_{5}^{\mathrm A_4}}
\newcommand{\Sing}{\operatorname{Sing}}
\newcommand{\Supp}{\operatorname{Supp}}
\newcommand{\g}{\operatorname{g}}

\newcommand{\Pic}{\operatorname{Pic}}
\newcommand{\Diff}{\operatorname{Diff}}
\newcommand{\Bs}{\operatorname{Bs}}
\newcommand{\ind}{\operatorname{Ind}}
\newcommand{\Cl}{\operatorname{Cl}}
\newcommand{\red}{\operatorname{red}}
\newcommand{\rk}{\operatorname{rk}}
\newcommand{\dif}{\operatorname{d}_{\Sha}}

\newcommand{\qq}{\mathbin{\sim_{\scriptscriptstyle{\QQ}}}}

\newcommand{\comment}[1]{}

\newcommand{\xref}[1]{{\rm \ref{#1}}}

\renewcommand\labelenumi{\rm (\arabic{enumi})}
\renewcommand\theenumi{\rm (\arabic{enumi})}
\renewcommand{\emptyset}{\varnothing}
\newtheorem{theorem}{Theorem}

\numberwithin{theorem}{section}
\numberwithin{equation}{theorem}

\newtheorem{mtheorem}[theorem]{} 
\newtheorem{stheorem}[equation]{}

\theoremstyle{definition}
\newtheorem{case}[theorem]{}
\newtheorem{scase}[equation]{}

\newtheorem*{remark*}{Remark}

\newcounter{NN}\numberwithin{NN}{section}
\renewcommand{\theNN}{\arabic{NN}${}^o$}
\def\nr{\refstepcounter{NN}{\theNN}}%

\begin{document}

\title{$\QQ$-Fano threefolds of index $7$}
\author{Yuri Prokhorov}
\thanks{This work is supported by the Russian Science Foundation under grant 14-50-00005.}

\address{Steklov Mathematical Institute, Russian Academy of Sciences
}
\email{prokhoro@mi.ras.ru}

\begin{abstract}
We show that, for a $\mathbb Q$-Fano threefold $X$ of Fano index 7, the inequality $\dim |-K_X| \ge 15$ 
implies that $X$ is isomorphic to one of the following varieties
$\mathbb P (1^2,2,3)$, 
$X_6 \subset \mathbb P (1,2^2,3,5)$ or $X_6 \subset \mathbb P (1,2,3^2,4)$.
\end{abstract}

\maketitle

\section{Introduction}
The purpose of this note is to demonstrate an application of the birational technique developed in
series of papers \cite{Prokhorov-2007-Qe}, \cite{Prokhorov2008a}, \cite{Prokhorov-2013-fano}, \cite{Prokhorov-Reid}
to the biregular classification of singular Fano threefolds.
Recall that a projective algebraic variety $X$
called \textit {$\QQ$-Fano} if it has only
terminal $\QQ$-factorial singularities
$\Pic (X)\simeq\mathbb Z$,
and the anticanonical divisor $-K_X$ is ample. The interest to study these varieties 
is justifyed by the fact that  
they are naturally appear as a result of the application of the minimal model program.
Recall definitions of Fano-Weil and $\QQ$-Fano indices:
\[
\begin{array}{ll}
\qW(X) &:=\max \{q \in \ZZ \mid \hbox {$-K_X \sim qA$, $A$ is Weil divisor}\}
\\[3pt]
\qQ(X) &:=\max \{q \in \ZZ \mid \hbox {$-K_X \qq qA$, $A$ is Weil divisor}\}.
\end{array}
\]
It is clear that $\qW (X)$ divides $\qQ(X)$ and these numbers can differ
only if the divisor class group $\Cl(X)$ has a torsion.
Below we assume that $X$ is three-dimensional.
It is well known (see \cite{Suzuki-2004}, \cite{Prokhorov2008a}) that the
index $\qQ(X)$ can take only  the following values
\addtocounter{theorem}{1}
\begin{equation}
\label{eq!Kaori}
\qQ(X) \in \{1, \dots, 9, 11, 13, 17, 19\}.
\end{equation}
$\QQ$-Fano threefolds  of $\QQ$-Fano indices $\qQ(X)\ge9$ studied more or less completely:
the standard arguments with using the orbifold Riemann-Roch formula \cite[\S 10]{Reid-YPG1987} and
Bogomolov-Miyaoka inequality (see \cite{Kawamata-1992bF}) allow to find all the numerical invariants.
Thus there are 10 numerical ``candidate varieties'', i.e. collections of numerical invariants (see \cite{GRD}, \cite[Prop. 3.6]{Prokhorov2008a}).
Among them, one case (with $\qQ(X)=10$) is not realized geometrically \cite{Prokhorov2008a} and the the remaining nine cases
examples \cite{Brown-Suzuki-2007j} are known. Furthermore,
in five cases, the corresponding $\QQ$-Fano threefolds are \textit {completely described}
\cite{Prokhorov2008a}, \cite{Prokhorov-2013-fano}.
The situation becomes more complicated for larger values of the index $\qQ(X)$. For example, 
for $\qQ(X)=8$ there are 10 ``numerical candidates'', five  of them are not realized geometrically
 \cite{Prokhorov-2013-fano}, in three cases examples are known  
\cite{Brown-Suzuki-2007j}, and
two of them are completely described \cite{Prokhorov-2013-fano}.
For  $\qQ(X)=7$ there are already 23 ``numerical candidates'' (see Table \ref{tabular}).
It is known that 7 of them do not occur (see \cite{Prokhorov-2007-Qe}, \cite{Prokhorov-2013-fano}),
in five cases examples are known \cite{Brown-Suzuki-2007j},
and two of them  are completely described \cite{Prokhorov-2013-fano}.
In this article, we  improve these results:

\addtocounter{theorem}{-1}
\begin{mtheorem}{\bf Theorem.} \label{mtheorem}
Let $X$ be a $\QQ$-Fano threefold with $\qQ (X)=7$.
If $\dim|-K_X |\ge15$, then
$X$ is isomorphic to one of the following varieties:
\[
\text{
$\PP (1^2,2,3)$, \quad
$X_6 \subset \PP (1,2^2,3,5)$ or $X_6 \subset \PP (1,2,3^2,4)$.}
\]
\end{mtheorem}

In fact, Theorem \ref{mtheorem} is a direct consequence of
the papers \cite{Prokhorov2008a}, \cite{Prokhorov-2013-fano}, and the following
technical result.

\begin{mtheorem}{\bf Theorem.}
\label{proposition-q=7main}
Let $X$ be a $\QQ$-Fano threefold of index $\qQ(X)=7$. Then, for $X$
the cases \xref{no-q=7-B23311} and \xref{no-q=7-B22510} from the Table \xref{tabular} do not occur.
In the case \xref{no-q=7-B2334}, the variety
$X$ can be embedded to $\PP (1,2,3^2,4)$ as a hypersurface of degree $6$.
\end{mtheorem}

\begin{scase}{\bf Remark.}
In the case \xref{no-q=7-B2334}, the  hypersurface $X=X_6 \subset \PP (1,2,3^2,4)$
has one of the following forms
\begin{eqnarray}\label{eqnarray-hypersurface-1}
y_2y_4+y_3^2+y_1^6+ y_3'\phi(y_1,y_2,y_3)&=&0,
\\
\label{eqnarray-hypersurface-2}
y_1^2y_4+y_3^2+y_2^3+ y_3'\phi(y_1,y_2,y_3)&=&0,
\end{eqnarray}
where $\phi$ is a polinomial of weighted degree $3$.
The index $4$ point 
is  a cyclic quotient $\frac 14 (1,1, -1)$ in the case  \eqref{eqnarray-hypersurface-1}
and a terminal singularity of type $\mathrm {cAx/4}$ in the case  \eqref{eqnarray-hypersurface-2} (see \S \ref{section-last}).
\end{scase}

\begin{scase}{\bf Remark.}
Thus, the question on the existence of $\QQ$-Fano threefolds with $\qQ(X)=7$ remains open in
the following nine cases:
\ref{no-q=7-B},
\ref{no-q=7-B369} -
\ref{no-q=7-B22258},
\ref{no-q=7-B2222512} -
\ref{no-q=7-B389}.

Note that in the cases
\ref{no-q=7-B22311} and
\ref{no-q=7-B22258}, the
expected threefold has codimension $4$ (see \cite[No. 41476, 41475]{GRD}).
It would be interesting to try to construct examples of using the methods of the work
\cite{Brown-Kerber-Reid-codim4}. It would be interesting also to study 
actions of finite groups on $\QQ$-Fano threefolds discussed here (cf. \cite{Prokhorov-planes}, \cite{PrzyjalkowskiShramov2016}).
\end{scase}

I thank the referee for numerous comments which 
helped me to improve the manuscript.

\section{Preliminaries}
Throughout this paper the ground field $\Bbbk$ is supposed to be 
algebraically closed of characteristic $0$.
We use the notation of the papers
\cite{Prokhorov-2007-Qe}, \cite{Prokhorov2008a}, \cite{Prokhorov-2013-fano}.
In particular, $\B(X)$ is the \textit{basket} of singularities of a terminal threefold $X$.

\begin{table}
\caption{\small $\QQ$-Fano threefolds with $\qQ(X)=7$}
\label{tabular}
\renewcommand{\arraystretch}{1.3}
\scalebox{0.8}{
\begin{tabular}{|l|l|l|l|l|l| |l|l|l|l|l|l|}
\hline
No.& $\B(X)$ &$A^3$&$\g(X)$&$\exists?$&ref&No.& $\B(X)$ &$A^3$&$\g(X)$&$\exists?$&ref
\\
\hline

\nr\label{no-q=7-B39} 
& $(3, 9)$
& $2/9$
&$38$
&$-$
& \cite{Prokhorov-2007-Qe}

&
\nr\label{no-q=7-B210} 
& $(2, 10)$
& $1/5$
&$34$
&$-$
& \cite{Prokhorov-2007-Qe}

\\

\nr\label{no-q=7-B23} 
& $(2, 3)$
& $1/6$
&$29$
&
$+!$ & \cite{Prokhorov2008a}

&
\nr\label{no-q=7-B2259} 
& $(2^2, 5, 9)$
& $7/45$
&$26$
& $-$
& \cite{Prokhorov-2013-fano}

\\
\nr\label{no-q=7-B2368} 
& $(2, 3, 6, 8)$
& $1/8$
&$21$
&$-$
& \cite{Prokhorov-2013-fano}

&
\nr\label{no-q=7-B228} 
& $(2^2, 8)$
& $1/8$
&$21$
&$-$
& \cite{Prokhorov-2013-fano}

\\
\nr\label{no-q=7-B2359} 
& $(2, 3, 5, 9)$
& $11/90$
&$20$
&$-$
& \cite{Prokhorov-2013-fano}

&
\nr\label{no-q=7-B23410} 
& $(2, 3, 4, 10)$
& $7/60$
&$19$
&$-$
& \cite{Prokhorov-2013-fano}

\\
\nr\label{no-q=7-B2225} 
& $(2^3, 5)$
& $1/10$
&$17$
&
$+!$ & \cite{Prokhorov-2013-fano}

&
\nr\label{no-q=7-B23311} 
& $(2, 3^2, 11)$
& $7/66$
&$17$
&$-$
&Th. \xref{proposition-q=7main}

\\
\nr\label{no-q=7-B22510} 
& $(2^2, 5, 10)$
& $1/10$
&$16$
&$-$
&Th. \xref{proposition-q=7main}

&
\nr\label{no-q=7-B2334} 
& $(2, 3^2, 4)$
& $1/12$
&$14$
&
$+!$ &Th. \xref{proposition-q=7main}

\\
\nr\label{no-q=7-B} 
& $(2^2, 3, 12)$
& $1/12$
&$13$
&$?$& 

&
\nr\label{no-q=7-B2235} 
& $(2^2, 3, 5)$
& $1/15$
&$11$
&
$+$ & \cite{Brown-Suzuki-2007j}

\\
\nr\label{no-q=7-B369} 
& $(3, 6, 9)$
& $1/18$
&$9$
&$?$& 

&
\nr\label{no-q=7-B22348} 
& $(2^2, 3, 4, 8)$
& $1/24$
&$6$
&$?$& 

\\
\nr\label{no-q=7-B2610} 
& $(2, 6, 10)$
& $1/30$
&$5$
&$?$& 

&
\nr\label{no-q=7-B22311} 
& $(2^2, 3, 11)$
& $1/33$
&$4$
&$?$& 

\\
\nr\label{no-q=7-B22258} 
& $(2^3, 5, 8)$
& $1/40$
&$3$
&
$?$& 

& 
\nr\label{no-q=7-B222345} 
& $(2^3, 3, 4, 5)$
& $1/60$
&$2$
&
$+$ & \cite{Brown-Suzuki-2007j}

\\
\nr\label{no-q=7-B2222512} 
& $(2^4, 5, 12)$
& $1/60$
&$1$
&$?$& 

&
\nr\label{no-q=7-B2313} 
& $(2, 3, 13)$
& $1/78$
&$1$
&$?$& 

\\
\nr\label{no-q=7-B389} 
& $(3, 8, 9)$
& $1/72$
&$1$
&$?$& 
&&&&&&
\\
\hline
\end{tabular}
}
\end{table}

\begin{mtheorem}{\bf Proposition (\cite[Lemma 6.1]{Prokhorov-2013-fano}, \cite{GRD}).}
Let $X$ be a $\QQ$-Fano threefold with $\qQ(X)=7$.
Then $\qW (X)=7$ and $X$ belongs to one of the numeric types in Table \xref{tabular},
where $\g (X):=\dim|-K_X |-1$ is the \emph{genus} of our Fano threefold $X$ and $A$ is a Weil divisor
on  $X$ such that $-K_X= \qW(X) A$.
\end{mtheorem}

We recall the construction used in the papers \cite{Prokhorov-2007-Qe},
\cite{Prokhorov2008a}, \cite{Prokhorov-2013-fano}, \cite{Prokhorov-Reid}.
\begin{case}
Let $X$ be a $\QQ$-Fano threefold.
For simplicity, we assume that the group $\Cl(X)$ is torsion free.
This holds in all our cases.
Everywhere throughout this paper by $A$ we denote
the positive generator of $\Cl(X)\simeq\ZZ$.
Then $-K_X=qA$, where $ q=\qQ(X)=\qW (X)$.

Consider a linear system $\MMM$ on $X$ without fixed components.
Let $c=\operatorname {ct} (X, \MMM)$ be the canonical threshold of the pair $(X, \MMM)$.
Consider a log crepant blowup $f: \tilde X \to X$ with respect to $K_X+c \MMM$
(see \cite{Alexeev-1994ge}).
Let $E$ be the exceptional divisor.
According to \cite{Alexeev-1994ge} one can choose $f$ so that
$\tilde X$ has only terminal $\QQ$-factorial singularities.
We can write
\begin{equation} \label{equation-1}
\begin{array}{lll}
K_{\tilde X} &\qq & f^*K_X+\alpha E,
\\[2pt]
\tilde\MMM &\qq & f^*\MMM- \beta E.
\end{array}
\end{equation}
where $\alpha \in \QQ_{>0}$, $\beta \in \QQ_{\ge0}$.
Then $c=\alpha/\beta$.
Assume that the log divisor $-(K_X+c \MMM)$ is ample.
Apply the log minimal model program with respect to $K_{\tilde X}+c \tilde\MMM$.
We obtain the following diagram
\begin{equation}\label{diagram-main}
\vcenter{
\xymatrix@C=19pt{
**[l] \tilde X=X_1\ar[d]_f\ar@/^19pt/@{-->}[rrrrrr]^{\chi}\ar@{-->}[r]&
\cdots\ar@{-->}[r]
&X_k\ar@{-->}[rr]_{\chi_k}\ar[rd]_{\pi_k}
&&X_{k+1}\ar@{-->}[r]\ar[ld]^{\mu_k}&\cdots\ar@{-->}[r]&**[r]\bar{X}\ar[d]^{\bar{f}}
\\
**[l] X&\cdots&&Y_k&&\cdots&**[r]\hat{X}
} 
}
\end{equation}
Here $\chi$ is a composition of $K_{\tilde X}+c \tilde\MMM$-log flips $\chi_k$,
all the varieties $X_k$
have only terminal $\QQ$-factorial singularities, 
$\uprho (X_k)=2$, and
$\bar{f}: \bar{X}=X_N \to \hat{X}$ is an extremal $K_{\bar{X}}$-negative Mori contraction.
In particular, $\rk \Cl(\hat{X})=1$.
In what follows, for the divisor (or linear system) $D$ on  $X$
by $\tilde D$ and $\bar D$ we denote
proper transforms of $D$
on  $\tilde X$ and $\bar{X}$ respectively.
If $|kA|\neq \emptyset$, we put $\MMM_k:=| kA|$
(is it possible that $\MMM_k$ has fixed components).
If $\dim \MMM_k=0$, then by $M_k$ we denote
a unique effective divisor $M_k \in \MMM_k$.
As in \eqref {equation-1}, we write
\begin{equation} \label{equation-12}
\bar\MMM_k\qq f^*\MMM_k- \beta_k E.
\end{equation}
\end{case}

\begin{stheorem}{\bf Lemma.} \label{lemma-flips-E}
All the maps $\chi_k: X_k \dashrightarrow X_{k+1}$ are
$-E_k$-flips \textup(flips with respect to the proper transform of $-E$\textup).
\end{stheorem}
\begin{proof}
By induction: on each step  the divisor $E_k$ is $\mu_{k-1}$-negative.
Since it is effective, it must be $\pi_k$-positive.
That means that $E_{k+1}$ is $\mu_{k}$-negative  by the definition of flips.
\end{proof}

\begin{case}
Assume that the contraction $\bar{f}$ is birational.
Then $\hat{X}$ is a $\QQ $-Fano threefold.
In this case, denote by $\bar F$ the
$\bar{f}$-exceptional divisor, by
$\tilde F \subset \tilde X$ its proper transform, $F:=f(\tilde F)$, and
$\hat{q}:=\qQ (\hat{X})$.
Let $\Theta$ be an ample Weil divisor on $\hat{X}$
generating $\Cl(\hat{X})/\operatorname {tors}$.
We can write
\[
\hat{E} \qq e\Theta,\qquad \hat\MMM_k \qq s_k\Theta,
\]
where $e\in \ZZ_{>0}$, $s_k \in \ZZ_{\ge0}$. If $\dim \MMM_k=0$ and $\bar M_k=\bar F$
(i.e. a unique element $M_k$ of the linear system $\bar\MMM_k$ is the $\bar{f}$-exceptional divisor), we  put
$s_k=0$.
\end{case}

\begin{case}
Assume that the contraction $\bar{f}$ is not birational.
In this case, $\Cl(\hat{X})$ has no torsion. Therefore, $\Cl(\hat{X})\simeq\ZZ$.
Denote by $\Theta$ the ample generator of $\Cl(\hat{X})$ and by
$\bar F$ a general geometric fiber.
Then $\bar F$ is either a smooth rational curve or a del Pezzo surface.
The image of the restriction map $\Cl(\bar{X}) \to \Pic (\bar F)$ is isomorphic to $\ZZ$.
Let $\Xi$ be its ample generator.
As above, we can write
\[
-K_{\bar{X}} |_{\bar F} =-K_{\bar F} \sim \hat{q}\Xi,\qquad \bar E |_{\bar F} \sim e \Xi,\qquad \bar\MMM_k |_{\bar F} \sim s_k \Xi,
\]
where $e\in \ZZ_{>0}$, $s_k \in \ZZ_{\ge0}$, and $\hat{q}\in \{1,2,3\}$.

If $\hat{X}$ is a curve, then $\hat{q}\le 3$ and $\hat{X}\simeq\PP^1$.
If $\hat{X}$ is a surface, then $\hat{q}\le 2$.
In this case, $\hat{X}$ can have only Du Val singularities of type
$\mathrm {A_n}$ \cite[Theorem 1.2.7]{Mori-Prokhorov-2008}.

\begin{stheorem}{\bf Lemma (see, e.g., \cite[Lemmas 3 \& 7]{Miyanishi-Zhang-1988}).}
\label{lemma-MZ}
Let  $\hat X$ be a del Pezzo surface whose singularities are at worst Du Val of type
$\mathrm{A_n}$. Assume that $\Cl(\hat X)\simeq \ZZ$. Then $\hat X$
is isomorphic to one of the following:
\[
\text{
$\PP^2$, $\PP (1,1,2)$, $\PP (1,2,3)$
or $\DP$,}
\]
where $\DP$ is a unique up to isomorphism
del Pezzo surface of degree $5$ whose singular locus consists of one point of type $\mathrm {A_4}$.
\end{stheorem}
\end{case}

Since  the group $\Cl(\bar X)$  has no torsion, the relations \eqref {equation-1} and \eqref {equation-12} give us
\begin{equation*} \label{equation-12-divisors}
k K_{\bar{X}}+q \bar\MMM_k \sim (k \alpha -q \beta_k) E.
\end{equation*}
From this we obtain
\begin{equation} \label{equation-main}
k \hat{q}=q s_k+(q \beta_k-k \alpha) e.
\end{equation}

\begin{mtheorem}{\bf Lemma (see \cite{Prokhorov2008a}).}
\label{corollary-qQ9}
Let $X$ ~ be a $\QQ$-Fano threefold
with $\qQ(X)\ge9$.
Then $\Cl(X)\simeq\mathbb Z$ and $\qQ(X)=\qW (X)$.
\end{mtheorem}

\begin{mtheorem}{\bf Lemma (see \cite[Lemma 4.2]{Prokhorov2008a}).}
\label{lemma-cthreshold}
Let $P \in X$ be a point of index $r>1$. Assume that
$\mathscr M \sim -mK_X$ near $P$, where $0<m<r$. Then $c \le 1/m$.
\end{mtheorem}

\begin{mtheorem}{\bf Theorem (\cite{Kawamata-1996}).}
\label{theorem-Kawamata-blowup}
Let $(Y \ni P)$ be a terminal quotient singularity of type
$\frac1r (1, a, r-a)$, let $ f \colon \tilde Y \to Y $ be a divisorial Mori contraction,
and let $E$ be the exceptional divisor.
Then $f(E)=P$, $f$ is a weighted blowup with weights $(1, a, r-a)$,
and the discrepancy of $E$ equals $a (E, X)=1/r$.
\end{mtheorem}

Note that in the case $\dim \hat{X}=3$, the group $\Cl(\hat{X})$
can have torsion elements. However, they can be easily controlled:

\begin{mtheorem}{\bf Lemma (see \cite[Lemma 4.12]{Prokhorov2008a}).}
\label{lemma-torsion-d}
Suppose that the map $\bar{f}$ is birational.
Write $F \sim dA$.
Then $\Cl(\hat{X})\simeq\mathbb Z \oplus \mathbb Z_{n}$, where $ d=ne$.
In particular, if $s_1=0$ \textup(i.e. $\dim \MMM_1=0$ and the divisor $\bar M_1$\enskip
is $\bar{f}$-exceptional\textup), then the group $\Cl(\hat{X})$
is torsion free and $e=1$.
\end{mtheorem}

\begin{case}
Suppose that the morphism $\bar{f}$ is birational.
Similar to \eqref {equation-1} and \eqref {equation-12} we can write
\[
K_{\bar{X}}\qq \bar{f}^*K_{\hat{X}}+b \bar F, \quad
\bar\MMM_k\qq \bar{f}^*\hat\MMM_k-\gamma_k \bar F, \quad
\bar E\qq \bar{f}^*\hat{E}-\delta \bar F.
\]
This gives us
\begin{equation*}
\begin{array}{lll}
s_kK_{\bar{X}}+\hat{q}\bar\MMM_k &\sim & (b s_k- \hat{q}\gamma_k) \bar F,
\\[2pt]
eK_{\bar{X}}+\hat{q}E &\sim & (be-\hat{q}\delta) \bar F.
\end{array}
\end{equation*}
Taking proper transforms of these relations to $X$, we obtain
\begin{equation} \label{equation-b-gamma-delta}
\begin{array}{lll}
-s_kq+k \hat{q}&=& (b s_k- \hat{q}\gamma_k) d,
\\[2pt]
-eq &=& (be-\hat{q}\delta) d.
\end{array}
\end{equation}
If $\Cl(\hat{X})\simeq\ZZ$, then $e=d$ by Lemma \ref{lemma-torsion-d}. Hence, 
\begin{equation} \label{equation-b-gamma-delta-1}
\begin{array}{lll}
be &=& \hat{q}\delta-q,
\\[2pt]
e \gamma_k &=& s_k \delta-k.
\end{array}
\end{equation}
\end{case}

\begin{mtheorem}{\bf Lemma.}
Let $X$ be a $\QQ$-Fano threefold of type \xref{no-q=7-B23311}, \xref{no-q=7-B22510} or
\xref{no-q=7-B2334}.
Then $\Cl(X)\simeq\ZZ$.
\end{mtheorem}
\begin{proof}
Suppose that $\Cl(X)$ contains an $n$-torsion element.
We may assume that $n$ is prime.
According to \cite[Prop. 2.9]{Prokhorov2008a} $n \le 7$ and 
for the basket $\B(X)=(r_1, \dots, r_l)$
the following condition holds: 
\[
\sum_{r_i \equiv 0 \mod n} r_i\ge16.
\]
In each case of \xref{no-q=7-B23311} -
\xref{no-q=7-B2334}, it is not satisfied.
\end{proof}

\begin{case}
We will always use, without additional reminder, the following lemma which
immediately follows from the  orbifold version of the Riemann-Roch theorem \cite[\S 10]{Reid-YPG1987}.

\begin{stheorem}{\bf Lemma.} \label{lemma-dimensions}
In cases \xref{no-q=7-B23311}, \xref{no-q=7-B22510} and
\xref{no-q=7-B2334} the dimension the linear systems $|kA|$ is 
given by  the following table:
\renewcommand{\arraystretch}{1.3}
\begin{center}
\rm
\scalebox{1}{
\setlongtables\renewcommand{\arraystretch}{1.3}
\begin{tabular}{|c|c|c|c|c|c|c|c|c|c|c|}
\hline
&&&\multicolumn{8}{c|}{$\dim|kA|$}
\\
\hhline{|~|~|~|--------}
&$\B$&$A^3$& $|A|$&$|2A|$&$|3A|$&$|4A|$&$|5A|$&$|6A|$&$|7A|$

\\
\hline
&&&&&&&&&
\\[-17pt]
\xref{no-q=7-B23311}&
$(2, 3, 3, 11)$
&$7/66$&$0$&$1$&$3$&$5$&$8$&$13$&$18$

\\
\xref{no-q=7-B22510}&
$(2, 2, 5, 10)$
&$1/10$&$0$&$1$&$2$&$5$&$8$&$12$&$17$

\\
\xref{no-q=7-B2334}&
$(2, 3, 3, 4)$
&$1/12$&$0$&$1$&$3$&$5$&$7$&$11$&$15$\\
\hline
\end{tabular}
}
\end{center}
\end{stheorem}
\end{case}

\begin{mtheorem}{\bf Lemma.} \label{lemma-Fano-dimA-ge3}
Let $X$ be a $\QQ$-Fano threefold such that $\qQ(X)=\qW (X)\ge3$, and
$\dim|A|\ge3$, where $-K_X \sim \qW (X) A$. Then $X$ is isomorphic to either $\PP^3$, a
quadric $Q \subset \PP^4 $ or a hypersurface $X_3 \subset \PP (1^4, 2)$.
\end{mtheorem}
\begin{proof}
Using a computer search (see \cite{GRD}, \cite{Suzuki-2004}, or 
\cite[Lemma 3.5]{Prokhorov2008a}),
we obtain $\qW (X) \le 4$. Furthermore,
$\g(X)\ge 21$ for $\qW (X)=3$ and $\g (X)\ge33$ for $\qW (X)=4$.
Now the result follows from \cite[Theorem 1.2]{Prokhorov-2013-fano}.
\end{proof}

\section{The case \ref{no-q=7-B23311}.}
In this section
we assume that
$X$ is a $\QQ$-Fano threefold of type \ref{no-q=7-B23311},
i.e. $-K_X=7A$, $A^3=7/66$, and $\B(X)=(2,3,3,11)$.
\begin{case}
Since $\dim|A|=0$ and $\dim|2A|=1$, the pencil $\MMM:=\MMM_2=| 2A|$ has no fixed components.
In a neighborhood of the index $11$ point we can write $A\sim -8K_X$ and $\MMM \sim -5K_X$.
By Lemma \ref{lemma-cthreshold} we have $c \le 1/5$.
Since $c=\alpha/\beta_2$, the following inequalities hold:
\begin{equation} \textstyle
\beta_2\ge5 \alpha,\qquad \beta_1\ge\frac 52 \alpha.
\end{equation}
The relation \eqref {equation-main} for $k=2$ has the form
\begin{equation} \label{equation-main-7}
2 \hat{q}=7 s_2+(7 \beta_2-2 \alpha) e\ge7 s_2+33e \alpha.
\end{equation}

\end{case}

\begin{stheorem}{\bf Lemma.} \label{lemma-B23311-first}
$f(E)$ is a point of index $11$.
\end{stheorem}
\begin{proof}
Assume the contrary. Then $\alpha\ge1/3$.
From \eqref {equation-main-7} we obtain $\hat{q}\ge6$.
So the contraction $\bar{f}$ is birational, $s_2>0$ (because $\MMM_2$
is a movable linear system), and $\hat{q}\ge9$.
In this case, $\Cl(\hat{X})\simeq\ZZ$ by Lemma \ref{corollary-qQ9}.
Below the table of dimensions of linear systems
\cite[Prop. 3.6]{Prokhorov2008a} is used. We obtain successively
$\dim|\Theta|\le 0$, $s_2\ge2$, and $\hat{q}\ge13$.
Then $\dim|2\Theta|\le 0$, $s_2\ge3$, and $\hat{q}\ge17$.
Continue successively: $\dim|4\Theta|\le 0$,
$s_2\ge5$, and $\hat{q}>19$. This contradicts \eqref {eq!Kaori}.
The assertion is proved.
\end{proof}

\begin{stheorem}{\bf Lemma.}
$\hat{X}\simeq\PP (1^3,2)$.
\end{stheorem}
\begin{proof}
Since $f(E)$ is a point of index $11$, we have $\alpha=1/11$ (by Theorem \ref{theorem-Kawamata-blowup}).
Then $A\sim -8K_X$ in a neighborhood of $f(E)$ and
so $\beta_1=8/11+m_1$, $\beta_2=5/11+m_2$,
where $m_i$ are non-negative integers.
The relation
\eqref{equation-main} for $k=1$ has the form
\begin{equation*}
\hat{q}=5e+7 (s_1+m_1e).
\end{equation*}
Hence, $\hat{q}\ge5$ and the contraction $\bar{f}$
is birational.
If $\hat{q}\neq 5$, then $\hat{q}\not\equiv 0 \mod 5$. Therefore, $s_1+m_1e>0$
and $\hat{q}\ge17$. Then
as in the proof of Lemma \ref{lemma-B23311-first}
we have $\Cl(\hat{X})\simeq\ZZ$, $|\Theta|=\emptyset$, $e\ge2$, $\hat{q}=17$, $s_1+2m_1=1$, and $s_1=1$.
In particular, $|\Theta|\neq \emptyset$.
The contradiction shows that $\hat{q}=5$. Then
$e=1$ and $s_1=0$.
By Lemma \ref{lemma-torsion-d} \ $\Cl(\hat{X})\simeq\ZZ$.
According to \cite[Th. 1.4 (vii)]{Prokhorov2008a} $\dim|\Theta|\le 2$.
Therefore, $s_k>1$ for $k>2$.
From \eqref {equation-main} for $k=2$ and $4$ we obtain
$s_2=1$ and $s_4=2$.
Hence, 
$\dim|\Theta|\ge1$, and
$\dim|2\Theta|\ge5$.
A computer search (see \cite{GRD} or \cite[Proof of Lemma 3.5]{Prokhorov2008a}) gives us that $\Theta^3=1/2$ and for
$\B(\hat{X})$
there are only two possibilities:
$\B(\hat{X})=(2)$ and $\B(\hat{X})=(2, 2, 3, 6)$.
According to \cite[Theorem 1.2 (v)]{Prokhorov-2013-fano} the latter possibility does not occur
and in the former one
$\hat{X}\simeq\PP (1^3,2)$.
\end{proof}

\begin{proof} [Proof of Theorem \xref{proposition-q=7main} in the case \xref{no-q=7-B23311}.]
The equations \eqref {equation-b-gamma-delta-1} give us
\begin{equation*}
b=5 \delta-7, \quad
\gamma_4=2 \delta-4.
\end{equation*}
In particular, $\delta\ge2$ and $b\ge3$. By Theorem \ref{theorem-Kawamata-blowup} \
$\bar{f}(\bar F)$ is a non-singular point $\hat{P} \in \hat{X}=\PP (1^3,2)$.
If $\delta=2$, then $\gamma_4=0$ and so $\bar\MMM_4=\bar{f}^*\hat\MMM_4=\bar{f}^*| 2\Theta|$.
But in this case $\dim \bar\MMM_4=6>\dim \tilde\MMM_4$. The contradiction shows that $\delta>2$ and $\gamma_4>0$.

Recall that $\hat\MMM_4 \subset|2\Theta|$.
Since $\dim \hat\MMM_4=5$, $\dim|2\Theta|=6$, and $\Bs|2\Theta|=\emptyset$, the linear subsystem $\hat\MMM_4\subset |2\Theta|$ 
is of codimension $1$ and consists of \textit {all} of elements passing through $\hat{P}$.
According to the main result of \cite{Kawakita-2001} the morphism
$\bar{f}$ is a weighted blowup of the point $\hat{P} \in \hat{X}$
with weights $(1, w_1, w_2)$, $\gcd(w_1,w_2)=1$ (for a suitable choice of local coordinates).

Now we introduce quasi-homogeneous coordinates $x_1, x_1', x_1'', x_2$ in $\hat{X}=\PP (1^3,2)$
so that $\hat{P}=(0: 0: 1: 0)$.
It is easy to see that sections $x_1'x_1''$, $x_1x_1''$, $x_2$ define  elements of $L_1, L_2, L_3 \in \hat\MMM_4\subset|2\Theta|$
which are transversal at $\hat{P}$.
This means that for any choice of local coordinates in $\hat{P}$ the multiplicity of a general
divisor $D \in \hat\MMM_4 \subset|2\Theta|$, given by the equation $\lambda_1x_1'x_1''+\lambda_2 x_1x_1''+\lambda_3x_2$,
at the point $\hat{P}$ with respect to the weight $(1, w_1, w_2)$ is equal to 1.
Thus, $\gamma_4=2 \delta-4=1$. The contradiction proves
Theorem \xref{proposition-q=7main}
in the case \xref{no-q=7-B23311}.
\end{proof}

\section{The case \ref{no-q=7-B22510}.}
\begin{case}
In this section
we assume that
$X$ is a $\QQ$-Fano threefold of type \ref{no-q=7-B22510},
i.e. $\Cl(X)=\ZZ \cdot A$, $-K_X=7A$, $A^3=1/10$, and $\B(X)=(2,2,5,10)$.
More precisely, from \cite[No. 41483]{GRD} we see
\begin{equation} \label{equationbB22510}
\textstyle
\B(X)=\left (2 \times \frac 12 (1,1,1), \quad \frac15 (1,2,3), \quad \frac {1}{10} (1,3,7) \right).
\end{equation}
\end{case}
\begin{case}
Take $\MMM:=| 3A|$.
In a neighborhood of the index $10$ point we can write $A\sim -3K_X$ and $\MMM \sim -9K_X$.
By Lemma \ref{lemma-cthreshold} the inequality $c \le 1/9$ holds.
Therefore,
\begin{equation} \label{equation-beta-B22510}
\beta_1\ge3 \alpha, \quad \beta_3\ge9 \alpha.
\end{equation}
The relation \eqref {equation-main} for $k=1$ has the form
\begin{equation} \label{equation-main-7-B22510}
\hat{q}=7 s_1+(7 \beta_1- \alpha) e\ge7 s_1+20e \alpha.
\end{equation}
\end{case}

\begin{mtheorem}{\bf Lemma.} \label{lemma7-B22510-index510}
$f (E)$ is a point of index of $5$ or $10$.
\end{mtheorem}
\begin{proof}
If $\alpha\ge1$, then 
from \eqref {equation-main-7-B22510} we obtain  $\hat{q}>19$, which contradicts \eqref {eq!Kaori}.
So (by Theorem \ref{theorem-Kawamata-blowup}) $f(E)$ is a point of index $r=2$, $5$ or $10$ and $\alpha=1/r$.
Suppose that $f(E)$ is a point of index $2$. Then $\alpha=1/2$ and
in a neighborhood of $f(E)$ we have $\MMM_1 \sim -K_X$.
Therefore, $\beta_1 \equiv \alpha \mod \ZZ$ and we can write
$\beta_1=1/2+m_1$,
where $m_1\ge1$ (see \eqref {equation-beta-B22510}).
From \eqref {equation-main-7-B22510} we obtain
\begin{equation*}
\hat{q}=3e+7 (s_1+m_1e).
\end{equation*}
Since $\hat{q}\neq 10$, we have $\hat{q}\ge13$. In particular, the contraction $\bar{f}$
is birational.
According to Lemma \ref{corollary-qQ9} we have $\Cl(\hat{X})\simeq\ZZ$.
If $\hat{q}=13$, then $s_1+m_1e=1$ and $e=2$.
Since $m_1\ge1$, this is impossible.
Thus, $\hat{q}\ge17$. Then $|\Theta|=\emptyset$ (see \cite[Prop 3.6.]{Prokhorov2008a}).
Therefore, $e\ge2$ and $s_1>0$ by Lemma \ref{lemma-torsion-d}.
But in this case, $\hat{q}>19$. This contradicts \eqref {eq!Kaori}.
The lemma is proved.
\end{proof}

\begin{mtheorem}{\bf Lemma.} \label{lemma-difficult-case}
Let $r$ be the index of the point $f(E)$.
There are the following possibilities:
\begin{enumerate}
\renewcommand\theenumi {\rm \themtheorem, \arabic {enumi})}
\renewcommand\labelenumi {\rm \arabic {enumi})}

\item
\label{lemma-difficult-case-2}
$r=10$, $\hat{X}$ is a surface, $e=1$, and linear systems
$\bar\MMM_1$, $\bar\MMM_2$, $\bar\MMM_3$
are vertical \textup(i.e. do not meet a general fiber\textup). There are two subcases:
\begin{enumerate}
\renewcommand\theenumii {{\rm \alph {enumii})}}
\renewcommand\labelenumii {{\rm \alph {enumii})}}

\item
\label{lemma-difficult-case-2-a}
$\hat{X}\simeq\PP (1,2,3)$,
\item
\label{lemma-difficult-case-2-b}
$\hat{X}$ is a 
del Pezzo surface of degree $5$ with a unique singularity which is  of type $\mathrm {A_4}$
\textup(see Lemma \xref{lemma-MZ}\textup).
\end{enumerate}

\item
\label{lemma-difficult-case-1}
$r=5$, $\hat{X}\simeq\PP (1,2,3,5)$, $e=1$, $s_1=0$, $s_3=3$.
\end{enumerate}

\end{mtheorem}

\begin{proof}
By Lemma \ref{lemma7-B22510-index510} either $r= 10$ or $r=5$.

Consider the case $r= 10$.
Then $\alpha=1/10$ and
$\beta_3=9/10+m_3$,
where $m_3$ is a non-negative integer.
The relation \eqref {equation-main} for $k=3$ is can be written as follows
\begin{equation} \label{equation-main-7-B22510-10}
3 \hat{q}=6e+7 (s_3+m_3e).
\end{equation}
In particular, $\hat{q}\ge2$.

Let $\hat{q}=2$. Then $e=1$ and $s_3=0$.
From \eqref {equation-main} for $k=1$ and $k=2$ we obtain $s_1=s_2=0$.
Since $\dim \bar\MMM_2=1$ and $\dim \bar\MMM_3=2$, the contraction $\bar{f}$ is not
birational and both linear systems $\bar\MMM_2$ and $\bar\MMM_3$ are
vertical.
If $\hat{X}\simeq\PP^1$, then $\bar\MMM_2=\bar{f}^*|\Theta|$
because $\dim \bar\MMM_2=1$. Similarly,
$\bar\MMM_3=\bar{f}^*| 2\Theta|\sim 2 \bar\MMM_2$,
a contradiction.
Therefore, $\hat{X}$ is a surface.
Since $\bar\MMM_1$ is a vertical divisor and $\dim|\bar\MMM_1|=0$,
we have $\hat{X} \not\simeq\PP^2$, $\PP (1,1,2)$.
We obtain the case \ref{lemma-difficult-case-2}.

Now let $\hat{q}>2$. Then $\hat{q}\ge4$ and so the contraction $\bar{f}$
is birational. In this case, $s_3>0$ (because $\dim \bar\MMM_3>0$).
From \eqref {equation-main-7-B22510-10} we obtain $\hat{q}\ge9$ and $s_3=3$.
Since $\dim \hat\MMM_3\ge2$ and $\dim|3\Theta|\le 1$ for $\hat{q}>11$
(see \cite[Proposition 3.6]{Prokhorov2008a}), 
there are exactly two possibilities: $\hat{q}=9$ and $\hat{q}=11$.
Moreover, $\dim|3\Theta|\ge2$.

If $\hat{q}=9$, then $\hat{X}\simeq X_6 \subset \PP (1,2,3,4,5)$ (see \cite[Proposition 3.6]{Prokhorov2008a}
and \cite[Th. 1.2]{Prokhorov-2013-fano}).
Then $\dim|7\Theta|=11$ and so $s_6\ge8$.
On the other hand, form \eqref {equation-main} for $k=6$, we have
\begin{equation*}
\textstyle
54=7 s_6+7 \beta_6- \frac 35,\qquad s_6 \le 7.
\end{equation*}
The contradiction shows that $\hat{q}\neq $ 9.

Let $\hat{q}=11$. Then $e=2$ and $\hat{X}\simeq\PP (1,2,3,5)$
(see \cite[Proposition 3.6, Th. 1.4]{Prokhorov2008a}).
Since $\hat\MMM_k \subset|s_k\Theta|$, from \eqref {equation-main} for $k=1$, $2$, $6$
we obtain $s_1=1$,
$s_2=2$, and
$s_6=8$. Moreover, comparing the dimensions of linear systems (Lemma \ref{lemma-dimensions} and
\cite[Proposition 3.6]{Prokhorov2008a}),
we obtain $\hat\MMM_2=| 2\Theta|$,
$\hat\MMM_3=| 3\Theta|$, and
$\hat\MMM_6=| 8\Theta|$.
The relations \eqref {equation-b-gamma-delta-1} can be written as
\begin{equation*}
\begin{array}{lll}
2b &=& 11 \delta-7,
\\[1pt]
2 \gamma_k &=& s_k \delta-k.
\end{array}
\end{equation*}
Hence, $\delta\ge1$ and $b\ge2$.
By Theorem \ref{theorem-Kawamata-blowup} \ $\bar{f}(\bar F)$ is a smooth point
and all the numbers $\delta$, $b$, $\gamma_k$ are integers.
If $\delta>1$, then $\gamma_1$, $\gamma_2$, $\gamma_3>0$.
Therefore, $\hat{P} \in \Bs|2\Theta|\cap \Bs|3\Theta|$.
But then $\hat{P}$ is a point of index $5$,
a contradiction. Therefore, $\delta=1$, $\gamma_1=\gamma_2=\gamma_3=0$, $\gamma_6=1$, and
\[
\hat{P} \in \Bs|8\Theta|=\{x_1=x_2=x_3x_5=0\},
\]
i.e. $\hat{P}$ is a point of index $5$ or $3$. This is again a contradiction.

Finally, consider the case $r=5$.
Then $\alpha=1/5$ and $\beta_3=9/5+m_1$,
where $m_3$ is a non-negative integer (see \eqref {equation-beta-B22510}).
The equality \eqref {equation-main} for $k=3$ takes the form
\begin{equation*}
3 \hat{q}=12e+7 (s_3+m_3e).
\end{equation*}
Taking the inequality $s_3>0$ into account as above we obtain $s_3=3$ and $\hat{q}=11$ or $19$.
In particular, $\dim|3\Theta|\ge2$.
If $\hat{q}=19$, then $\dim|3\Theta|=0$ \cite[Proposition 3.6]{Prokhorov2008a},
a contradiction.
Thus, $\hat{q}=11$.
Then, according to
\cite[Proposition 3.6, Theorem 1.4]{Prokhorov2008a}
we have $\hat{X}\simeq\PP (1,2,3,5)$.
This is the case \ref{lemma-difficult-case-1}.
\end{proof}

\begin{stheorem}{\bf Corollary.} \label{corollary-baskets-1}
We have
\[
\B(\tilde X) =
\begin{cases}
(\underline 2, \underline 2, \underline 5,3,7) & \text {in the case \xref{lemma-difficult-case-2}},
\\
(\underline 2, \underline 2,2,3, \underline {10}) & \text {in the case \xref{lemma-difficult-case-1}}.
\end{cases}
\]
where underlining $\underline {*}$ allocates the points contained in $\tilde X \setminus E$.
\end{stheorem}
\begin{proof}
Follows from \eqref {equationbB22510} and the fact that
$f$ is a weighted blowup of the point of index $r=10$ (respectively, $5$) with weights $\frac {1}{10} (1,3,7)$
(resp. $\frac 15 (1,2,3)$).
\end{proof}

\begin{stheorem}{\bf Lemma.} \label{lemma-Bs-K}
The divisor $-K_{\tilde X}$ is ample and the base locus of
$|-K_{\tilde X}|$ does not contain curves in $\tilde X \setminus E$.
\end{stheorem}
\begin{proof}
The variety $\tilde X$ is  FT type (Fano type) \cite[\S 2]{Prokhorov-Shokurov-2009}.
In particular, the ampleness property of $-K_{\tilde X}$ is equivalent to the fact that $-K_{\tilde X}$
has a positive intersection number with any curve $C \subset \tilde X$.
Suppose that $K_{\tilde X} \cdot C\ge0$.
If $C \subset E$, then $E \cdot C<0$, $ f^*K_X \cdot C=0$, and $K_{\tilde X} \cdot C=f^*K_X \cdot C+\alpha E \cdot C<0$.
So we may assume that $C \not\subset E$.

Conditions of Corollary 6.3 of  \cite{Prokhorov-2007-Qe} are are satisfied for our variety  $X$ in this case.
Therefore, the linear system $\MMM_7:=| -K_X|$ has only isolated base points.
We can write $K_{\tilde X}+\tilde\MMM_7+\lambda E\sim f^*(K_X+\MMM_7) \sim 0$,
where $\lambda=\beta_7- \alpha\ge0$.
Therefore, $-K_{\tilde X} \sim \tilde\MMM_7+\lambda E$. This proves the second assertion.
Since $C \not\subset E$, we have $(\tilde\MMM_7+\lambda E) \cdot C \le 0$. Since the
linear system $\tilde\MMM_7$
has only isolated base points outside of $E$, we have $\tilde\MMM_7 \cdot C=E \cdot C=0$.
But this is impossible because $\tilde\MMM_7$ and $E$ generate the group $\Cl(\tilde X)$.
\end{proof}

Now we  analyze the middle part of the diagram \eqref {diagram-main} in details.
Recall that $\chi$ is a composition of log flips $\chi_k: X_k \dashrightarrow X_{k+1}$.

\begin{stheorem}{\bf Corollary.} \label{lemma-flips}
All the maps $\chi_k: X_k \dashrightarrow X_{k+1}$ are
flips \textup(with respect to the canonical divisor\textup).
\end{stheorem}
\begin{proof}
Similar to the proof of Lemma \ref{lemma-flips-E}.
\end{proof}

\begin{mtheorem}{\bf Lemma.} \label{lemma-difficult-case-P (1,2,3,5)}
In the case \xref{lemma-difficult-case-1}, the contraction $\bar{f}$ is a
weighted blowup of a non-singular point $\hat{P} \in \PP (1,2,3,5)$
with weights $(1,1,3)$.
\end{mtheorem}

\begin{proof}
The relations \eqref {equation-b-gamma-delta-1} have the form
\begin{equation} \label{equation-b-gamma-delta-1-l}
\begin{array}{lll}
\gamma_3 &=& 3 \delta-3,
\\[1pt]
b &=& 11 \delta-7.
\end{array}
\end{equation}
Hence, $\delta\ge1$ and $b\ge4$.
Therefore, $\bar{f}(\bar F)$ is a smooth point
and all the numbers $\delta$, $b$, $\gamma_3$ are integers.
According to  \cite{Kawakita-2001} the morphism
$\bar{f}$ is a weighted blowup of the point $\hat{P} \in \hat{X}$
with weights $(1, w_1, w_2)$, where $\gcd (w_1, w_2)=1$ (for a suitable choice of local coordinates).
In particular, $\bar{X}$ has only (terminal) quotient singularities and
\begin{equation} \label{corollary-baskets-2}
\B(\bar{X})=(2,3,5, w_1, w_2).
\end{equation}
Recall that \textit {Shokurov's difficulty} $\dif(V)$ of a variety $V$ with terminal singularities
is
defined as the number of exceptional divisors on $V$ with the discrepancy $<1$ \cite[Definition 2.15]{Shokurov-1985}.
It is known that this number is well-defined and finite.
Moreover, in the three-dimensional case, it is strictly decreasing under flips \cite[Corollary 2.16]{Shokurov-1985}.
If $V \ni P$ is a terminal cyclic quotient of index $r$, then $\dif(V \ni P)=r-1$.

We claim that $\dif(\tilde X)= 14$. If the singularities of
$\tilde X$ are cyclic quotients,
this follows directly from the above stated and Corollary  \ref{corollary-baskets-1}.
Otherwise, $\tilde X$ has a unique non-quotient singularity  $\tilde P_2 \in \tilde X$
with $\B(\tilde X, \tilde P_2)=(2,2)$ (see Corollary \ref{corollary-baskets-1}).
Since the singularities of
$\bar{X}$ are cyclic quotients, the point $\tilde P_2$ should lie on a 
flipping curve. According to the classification flips \cite[Th. 2.2]{Kollar-Mori-1992}
$\tilde P_2$ is of type $\mathrm {cA/2}$.
Further, according to \cite[Th. 6.1, Rem. 6.4B]{Reid-YPG1987}, the singularity $\tilde P_2 \in \tilde X$
can be defined locally in the form 
\[
\{x_1x_2+\phi (x_3^2, x_4)=0\}/\mumu_2 \subset \CC^4/\mumu_2 (1,1,1,0),
\]
where $\operatorname {mult}_0 \phi (0, x_4)=2$. In this case, all the divisors
with discrepancy of $1/2$ are obtained as exceptional divisors of two explicitly  described
weighted blowups \cite{Kawamata-1992-e-app}.
This proves that $\dif(\tilde X)= 14$.
Thus, from the \eqref {corollary-baskets-2} we obtain
\[
\dif(\tilde X)=14\ge
\dif(\bar{X})=5+w_1+w_2=5+b=11 \delta-2.
\]
Hence, $\delta=1$, $b=4$, and $\gamma_3=0$.
Since $\gcd (w_1, w_2)=1$, we have (up to a permutation)
a unique possibility  $(w_1, w_2)=(1,3)$.
\end{proof}

\begin{stheorem}{\bf Lemma.} \label{lemma-cb}
Let $V$ be a threefold with
terminal singularities and let $\varphi: V \to W$ be a
contraction to a surface such that $-K_V$ is $\varphi$-ample and
all the fibers are of dimension $1$
\textup(i.e. $\varphi$ is a $\QQ$-conic bundle \textup).
Let $l$ be a general geometric fiber.
Assume that there is a $\varphi$-ample $\QQ$-Cartier Weil divisor $D$
such that $D \cdot l=1$.
Then the following assertions hold.
\begin{enumerate} [leftmargin=20pt]
\item
All fibers of $\varphi$ is irreducible.
\item
If the point $w \in W$ is smooth, then $\varphi$ is smooth morphism over $w$.
\item
If the point $w \in W$ is singular, it is Du Val of type $\mathrm {A_{n-1}}$
and $V$ has exactly two singular points on the fiber $\varphi^{- 1} (w)$.
These points are terminal quotients
of type $\frac 1n (1, a, -a)$ and $\frac 1n (-1, a, -a)$, $\gcd (n, a)=1$.
The support of the divisor $D$
contains at least one of these points.
\end{enumerate}
\end{stheorem}

\begin{proof}
Fix a point $w \in W$ and let $C:=\varphi^{- 1} (w)_{\mathrm {red}}$.
We can replace $W$ with a small analytic neighborhood of $w$.
Apply to $V$ over $(W \ni w)$ the analytic minimal model program:\
$V \dashrightarrow V'$ (see \cite[\S 4]{Nakayama1987}).
By our assumption $D \cdot l=1$, all the neighboring fibers are irreducible.
Therefore, $V \dashrightarrow V'$ is a sequence of flips.
On the other hand, since the general fiber $\varphi$ is a rational curve, the divisor $K_{V'}$ cannot be
nef over $W$.
Therefore, $V \dashrightarrow V'$ is an isomorphism.
Then $\uprho_{\mathrm {an}} (V/W)=1$ and the fiber over $w$ is irreducible.

Assume that the point $w \in W$ is smooth.
Let $P_1, \dots, P_l$ be all the singular points on $C$, and let $r_1, \dots, r_l $ be their indices.
Then the divisor $r_1 \cdots r_lD$ is Cartier in a neighborhood of $C$.
In particular, the intersection number $r_1 \cdots r_lD \cdot C$ is integral.
Suppose that at least one of the points $P_1, \dots, P_l$ is not Gorenstein.
Then according to \cite[Corollary 2.7.4, Lemma 2.8]{Mori-Prokhorov-2008}
all the numbers $r_1, \dots, r_l$ pairwise relatively prime and $-K_V \cdot C=1/(r_1 \cdots r_l)$
(because the base $W$ is smooth).
But then $r_1 \cdots r_lD \cdot C=1/2$.
The contradiction shows that $K_V$ is Cartier in a neighborhood of $C$
and the same is true for $D$.
Consider the scheme fiber $Z:=\varphi^{- 1} (w)$. As the
morphism is flat in a neighborhood of $w$, we have $D \cdot Z=1$.
Therefore, the fiber of $Z$ over $w$ is reduced.
According to the adjunction formula $Z=C$ is a smooth rational curve.
So morphism $\varphi$ smooth in this case.

Now let $w \in W$ be a singular point.
According to \cite[Theorem 1.2.7]{Mori-Prokhorov-2008} it is of type $\mathrm {A_{n-1}}$.
Therefore, $W$ is a  quotient of an open disc $0 \ni W'\subset \CC^2$ by
cyclic group $\mumu_n$ of order $n$, acting freely outside the origin.
Consider base change (see, e.g., \cite[2.4]{Mori-Prokhorov-2008})
\[
\xymatrix@R=10pt {
V'\ar [d]_{\varphi'} \ar [r] & V \ar [d]^{\varphi}
\\
W'\ar [r] & W
}
\]
Here $V'\to V$ is the quotient morphism by $\mumu_n$, where $\mumu_n$ acts on
$V'$ freely in codimension 2. The preimage $D' \subset V'$ of the section $D$
is also a section. As above we see that the morphism $\varphi': V' \to W'$ is
smooth. Then $\varphi$ is  toroidal
(see \cite[1.2.1]{Mori-Prokhorov-2008})
The action $\mumu_n$ on
$V'$ and the local structure of $V$ of the fiber near $w$ are completely described in
\cite[1.1.1]{Mori-Prokhorov-2008}.
In particular, $V$ has exactly two cyclic quotients
of type $\frac 1n (1, a, -a)$ and $\frac 1n (-1, a, -a)$ for some $a$
such that $\gcd (n, a)=1$.
The divisor $D$ cannot be  Cartier near $\varphi^{- 1} (w)$
(otherwise $\varphi$ is a smooth morphism). Therefore,
$D$ passes through at least one of the singular points.
\end{proof}

\begin{stheorem}{\bf Lemma.} \label{corollary-baskets-bar}
We have
\[
\B(\bar{X}) =
\begin{cases}
(2, \underline 2,3, \underline 3) & \text {in the case \xref{lemma-difficult-case-2-a}},
\\
(5, \underline 5) & \text {in the case \xref{lemma-difficult-case-2-b}},
\\
(\underline 2, \underline 3, \underline 5,3) & \text {in the case \xref{lemma-difficult-case-1}}.
\end{cases}
\]
where underlining $\underline {*}$ allocates the points lying on $\bar E$
\textup(but it is possible that some of non-underlining points also lie on $\bar E$\textup).
\end{stheorem}
\begin{proof}
In the case \xref{lemma-difficult-case-1},
the assertion follows from Lemma \ref{lemma-difficult-case-P (1,2,3,5)}.
Consider the case \xref{lemma-difficult-case-2}.
Since $e=1$ and $\uprho (\bar{X}/\hat{X})=1$, the
divisor $\bar E$ is a relatively ample birational section of the  morphism $\bar{f}$.
Put $U:=\hat{X} \setminus \Sing(\hat{X})$ and $V:=\bar{f}^{- 1} (U)$.
By Lemma \xref{lemma-cb} \ $V$ is smooth.
Now let $Q \in \hat{X}$ be a singular point of type $\mathrm {A_{n-1}}$ ($n=2$, $3$, $5$).
Again by Lemma \xref{lemma-cb} \ $V$ has exactly two cyclic quotients of
index $n$ and the divisor $\bar E$ passes through at least one of these points.
\end{proof}

\begin{stheorem}{\bf Corollary.}
The map $\chi$ is not an isomorphism.
\end{stheorem}
\begin{proof}
It is a consequence of Corollary \ref{corollary-baskets-1} and Lemma \ref{corollary-baskets-bar}.
\end{proof}

\begin{stheorem}{\bf Lemma.} 
In the case \xref{lemma-difficult-case-1}, the
linear system $|2\Theta|$ is the proper transform
linear system $|2A|$ \textup(in particular, $s_2=2$, and $\gamma_2=0$\textup).
\end{stheorem}
\begin{proof}
Consider the linear system $\hat\LLL:=| 2\Theta|$.
It is clear that $K_{\hat{X}}+\hat\LLL+9 \hat{E} \sim 0$.
Hence,
$K_{\bar{X}}+\bar\LLL+9 \bar E+a \bar F \sim 0$,
where $a\ge9 \delta-b=5$. Taking the proper transform of this relation to $X$, we obtain
$K_{X}+\LLL+aA\sim 0$ and $\LLL \sim (7-a) A$.
Since $s_1=0$, we have $\LLL \not\sim A$.
Since $\dim \LLL=1$, we have $\LLL \sim 2A$, $a=5 $, and $\LLL=\MMM_2$.
So, $s_2 =2$ and $\hat\MMM_2=| 2\Theta|$ as $\dim|2\Theta|=\dim \MMM_2$.
From \eqref {equation-b-gamma-delta-1-l} we obtain $\gamma_2=0$.
\end{proof}

\begin{stheorem}{\bf Corollary.}
The linear system $\bar\MMM_2$ is nef.
\end{stheorem}
\begin{proof}
Since $\gamma_2=0$, we have $\bar\MMM_2 \sim \bar{f}^*(2\Theta)$.
\end{proof}

\begin{stheorem}{\bf Lemma.}
All the maps $\chi_k: X_k \dashrightarrow X_{k+1}$ are
$\bar\MMM_2$-flips.
\end{stheorem}
\begin{proof}
Similar to the proof of Lemma \ref{lemma-flips-E} but using the descending induction on  $k$.
\end{proof}

Denote by $\ind (\chi)$ the indeterminacy locus of a map $\chi$.

\begin{stheorem}{\bf Corollary.} \label{corollary-BsM2}
$\ind (\chi) \subset \Bs|\tilde\MMM_2 |$.
\end{stheorem}

\begin{stheorem}{\bf Lemma.} \label{lemma-bs-2A}
The set $\Bs|2A|$ is irreducible.
\end{stheorem}

\begin{proof}
Let $M_1 \in | A|$ be the \textup(only \textup) effective
divisor and let $M_2 \in \MMM_2$ be a general member.
Consider the $1$-cycle $\Gamma:=M_1 \cap M_2$.
Since $\dim \MMM_2=1$, we have $\Supp (\Bs \MMM_2)=\Supp (\Gamma)$.
Write $\Gamma=\sum \gamma_i \Gamma_i$.
As
$ H:=10 A$ is an ample Cartier divisor,
we have
\[
2=20 A^3=H \cdot \Gamma=\sum \gamma_i H \cdot \Gamma_i\ge\sum \gamma_i.
\]
Therefore, $\Gamma$ has at most two irreducible components.
If $\Gamma$ has exactly two irreducible components, then
the inequality above becomes an equality and so $\gamma_1=\gamma_2=1$,
$\Gamma=\Gamma_1+\Gamma_2$, where $\Gamma_i \cdot A=1/10$.
In this case, the curve $\Gamma=M_1 \cap M_2$ is generically reduced. Therefore, the surface $M_2$
is smooth at the general points of curves $\Gamma_i$. By Bertini's theorem
$M_2$ is smooth outside of $\Gamma$. Therefore, $M_2$ is a normal surface.
Since the threefold $X$ is smooth in codimension 2, the usual
adjunction formula holds for $M_2$:
\[
K_{M_2}=(K_X+M_2) |_{M_2}=-5A|_{M_2}.
\]
Next we apply the adjunction formula for singular varieties
(see \cite[ch 16.]{Utah}) to the components of $\Gamma$:
\[
K_{\Gamma_1}+\Diff_{\Gamma_1} (\Gamma_2)=(K_{M_2}+\Gamma) |_{\Gamma_1} =
-4A|_{\Gamma_1},
\]
where $\Diff_{\Gamma_1} (\Gamma_2)$ is the \textit {different}, an effective divisor,
whose support is contained in $\Sing(M_2) \cup (\Gamma_1 \cap \Gamma_2)$.
Here
\[
\deg \Diff_{\Gamma_1} (\Gamma_2) =-\textstyle \frac 4 {10}-\deg K_{\Gamma_1}.
\]
Therefore, $\deg K_{\Gamma_1}<0$, $\Gamma_1\simeq\PP^1$
and $\deg \Diff_{\Gamma_1} (\Gamma_2)=8/5$.
We write $\Diff_{\Gamma_1} (\Gamma_2)=\sum \alpha_i P_i$, where $\alpha_i \in \QQ_{>0}$.
By the inversion of adjunction \cite[Theorem 17.6]{Utah} the
inequality $\alpha_i<1$ holds if and only if the pair
$(M_2, \Gamma)$ log terminal at $P_i$. In this case
$P_i \notin \Gamma_1 \cap \Gamma_2$ and
the singularity 
$M_2 \ni P_i$ is analytically isomorphic to the quotient $\frac 1n (1, q)$, $\gcd (n, q)=1$
(see \cite[Proposition 16.6]{Utah}).
On the other hand, the curve $\Gamma$ is an intersection of ample divisors and therefore connected.
Thus, we may assume that $\alpha_1\ge1$.
But then the only possibility is
\[\textstyle
\Diff_{\Gamma_1} (\Gamma_2)=\frac {11}{10} P_1+\frac {1}{2} P_2.
\]
In particular, $\Gamma_1 \cap \Gamma_2=\{P_1\}$.
On the other hand, $\Gamma$ contains  points $P'$, $P''$ indices $5$ and $10$ (because the divisors $M_1$ and $M_2$
are not Cartier at these points).
Since $P_2 \in M_2$ is a singularity of type $\frac12 (1,1)$
(see \cite[Proposition 16.6]{Utah}), we have $P_2 \notin \{P', \, P''\}$.
We may assume that $P'= P_1$ and $P''\notin \Gamma_1$.
By symmetry, $P'' \notin \Gamma_2$, a contradiction.
\end{proof}

\begin{stheorem}{\bf Lemma.} \label{lemma-flips-iso}
Let $C:=\Gamma_{\red}$ and let
$\tilde C \subset \tilde X$
be the proper transform of $C$.
Then
$\tilde X \setminus (E \cup \tilde C)\simeq\bar{X} \setminus \bar E$.
\end{stheorem}

\begin{proof}
According to Corollary \ref{corollary-BsM2} we have
$\tilde C \supset \ind (\chi)$.
Since $\chi$ is not an isomorphism, we have $\tilde C=\ind (\chi)$.
By Lemma \ref{lemma-flips-E} \ $\ind (\chi^{- 1}) \subset \bar E$.
\end{proof}

\begin{stheorem}{\bf Lemma.} \label{lemma-2-points}
The curve $\tilde C$
contains at least two non-Gorenstein points $P', \, P''\in \tilde X \setminus E$.
\end{stheorem}

\begin{proof}
In the case \xref{lemma-difficult-case-2}, the curve $\tilde C$
must contain the point of index $5$ (because this is holds for $C$).
On the other hand, in this case the set $\B(\tilde X \setminus E)$
contains two points of index $2$ and $\B(\bar{X} \setminus \bar E)$ contains 
at most one such a point (see Corollary
\ref{corollary-baskets-1} and Lemma \ref{corollary-baskets-bar}).
By Lemma \ref{lemma-flips-iso}
at least one index $2$ point lies on $\tilde C$.

Similarly, in the case \xref{lemma-difficult-case-1}, the set
$\tilde X \setminus E$
contains points of index $2$ and $10$, and the set $\bar{X} \setminus \bar E$
does not contain such points.
\end{proof}

\begin{proof} [Proof of Theorem \xref{proposition-q=7main} in the case \xref{no-q=7-B22510}.]
Let $D \in | -K_{\tilde X}|$ be a general anticanonical divisor.
By Lemma \ref{lemma-Bs-K}
$D \not\supset \tilde C$.
On the other hand, $D$ passes through all the non-Gorenstein points of $\tilde X$.
So the set $\Supp (D) \cap \tilde C$ is not connected.
This contradicts Shokurov's connectedness theorem  \cite[Th. 17.4]{Utah}
(as well as the classification of three-dimensional flips \cite[Th. 2.2]{Kollar-Mori-1992}).
Theorem \xref{proposition-q=7main}
is proved in the case \xref{no-q=7-B22510}.
\end{proof}

\section{The case \ref{no-q=7-B2334}}
\label{section-last}
In this section
we assume that
$X$ is a $\QQ$-Fano threefold of type \ref{no-q=7-B2334},
i.e. $\qW(X)=7$ and $\B(X)=(2,3,3,4)$.

\begin{case}{\bf Remark.}
For a point $P_4 \in X$ of index $4$ there are two possibilities:
\begin{enumerate}[leftmargin=20pt]
\renewcommand\theenumi{{\rm\alph{enumi})}}
\renewcommand\labelenumi{{\rm\alph{enumi})}}
\item
$P_4 \in X$ is cyclic quotient of type $\frac14 (1,1, -1)$,
in this case $X$  also has a point of type $\frac12 (1,1,1)$;
\item
$P_4 \in X$ is singularity of type $\mathrm {cAx/4}$ \cite[6.1 (2)]{Reid-YPG1987} and
$X$ has no points of index $2$.
\end{enumerate}
In both cases, the discrepancy of any divisorial blowup
$f: (\tilde X, E) \to (X, P_4)$ of $P_4$ in the Mori category
is equal to $1/4$ (see \cite{Kawamata-1996}, \cite{Kawakita2005}).
\end{case}

\begin{mtheorem}{\bf Proposition.} \label{proposition-no-q=7-B2334-c}
Let $X$ be a threefold of type \xref{no-q=7-B2334}.
Then the pair $(X, | 3A|)$
has only canonical singularities.
\end{mtheorem}

\begin{proof}
Assume that the pair $(X, \MMM=| 3A|)$ is not canonical.
Apply the construction \eqref {diagram-main}. Then $c=\alpha/\beta_3<1$.
Therefore,
\[
\beta_3>\alpha,\qquad 3 \beta_1>\alpha.
\]
The relation \eqref {equation-main} for $k=1$ and $3$ has the form
\begin{eqnarray} \label{equation-main-7-B2334-1a}
\hat{q}&=7 s_1+(7 \beta_1- \alpha) e &=7s_1+4 \beta_1e+(3 \beta_1- \alpha) e,
\\[1pt]
\label{equation-main-7-B2334-3a}
3 \hat{q}&=7 s_3+(7 \beta_3-3 \alpha) e &=7s_3+4 \beta_3e+3 (\beta_3- \alpha) e.
\end{eqnarray}
Suppose that $f$ is a blowup a curve or a Gorenstein point. Then
$\alpha$, $\beta_1$, $\beta_3$ are positive integers, $\beta_3\ge2$, and
$\hat{q}\ge7s_1+4 \beta_1+1\ge5$. Therefore, in this case the contraction $\bar{f}$ is birational.
If $s_1 \neq 0$, then $\hat{q}\ge13$.
If $s_1=0$, then $\Cl(X)$ is torsion free by Lemma \ref{lemma-torsion-d}.
In both cases, $\Cl(\hat{X})\simeq\ZZ$. By Lemma \ref{lemma-Fano-dimA-ge3} \ $\dim|\Theta|\le 2$.
Since $\dim|3A|=3$, we obtain successively
$s_3\ge2$,\
$3 \hat{q}\ge7s_3+11e\ge25$, \ $\hat{q}\ge9$, \ $s_3\ge4$, \ $\hat{q}\ge13$ \
$s_3\ge5$, \ $\hat{q}\ge17$, \ $s_3\ge9$, and, finally,
$\hat{q}>19$. This contradicts \eqref {eq!Kaori}.
Hence $f(E)$ is a point of index $r=2$, $3$ or $4$, and $\alpha=1/r$ or $\alpha=2/3$
(see \cite[Lemma 2.6]{Prokhorov-2013-fano}).

Consider the case where $f(E)$ is a point of index $4$.
Then $\alpha=1/4$, $A\sim -3K_X$ in a neighborhood of $f(E)$, and
so $\beta_1=3/4+m_1$, $\beta_3=1/4+m_3$,
where $m_i$ are non-negative integers, $m_3>0$.
The relations \eqref {equation-main-7-B2334-1a} and \eqref {equation-main-7-B2334-3a}
take the form
\begin{eqnarray*}
\hat{q}&=& 5e+7 (s_1+em_1),
\\[1pt]
3 \hat{q}&=& e+7 (s_3+em_3).
\end{eqnarray*}
We obtain a unique solution:
$\hat{q}=5$, $e=s_3=1$, $s_1=0$.
Moreover, $\Cl(\hat{X})\simeq\ZZ$. Therefore, $\dim|\Theta|\ge3$.
This contradicts \cite[Th. 1.4 (vii)]{Prokhorov2008a}.

Consider the case where $f(E)$ is a point of index $3$ and $\alpha=1/3$.
As above we obtain a unique solution:
$\hat{q}=2$, $e=1$, $s_1=s_3=0$ (see Lemma \ref{lemma-torsion-d}).
Since $\dim \MMM_3>0$ and $s_3=0$, the contraction $\bar{f}$ is not birational and the linear systems 
$\bar\MMM_1$, $\bar\MMM_3$ are
vertical.
Since $\dim \bar\MMM_1=0$, we have $\hat{X} \not\simeq\PP^2$, $\PP (1,1,2)$.
Since $\dim \bar\MMM_3=3$, the variety
$\hat{X}$ cannot be a surface according to \cite[Lemma 2.8]{Prokhorov-2013-fano}.
Therefore, $\hat{X}\simeq\PP^1$.
But then $\bar\MMM_3=\bar{f}^*|\OOO_{\PP^1} (3)|$ and
$\bar\MMM_1 \sim \frac 13 \bar\MMM_3 \sim f^*\OOO_{\PP^1} (1)$.
So,
$\bar\MMM_1=\bar{f}^*|\OOO_{\PP^1} (1)|$ is a movable linear system,
a contradiction.

Consider the case where $f(E)$ is a point of index $2$. Recall that
$\dim|3A|=3$. When for $\hat{q}\ge9$ the inequality $\dim|3\Theta|\le 2$ holds
(see \cite[Prop. 3.6]{Prokhorov2008a}).
Then $s_3\ge4$.
Using this, as above, we obtain a unique solution:
$\hat{q}=3$, $e=1$, $s_1=s_3=0$.
Since $s_3=0$ and $\hat{q}=3$, we have $\hat{X}\simeq\PP^1$.
But then, as above, $\bar\MMM_3=\bar{f}^*|\OOO_{\PP^1} (3)|$ and
$\dim \bar\MMM_1=1$,
a contradiction.

It remains to consider the case where $f(E)=P_3$ is point of index $3$ and $\alpha=2/3$.
Then a general element of the linear system $\MMM_3$ passes through 
$P_3$ and is a Cartier divisor at $P_3$.
As above we obtain the following possibilities:
$(\hat{q}, e, s_3)=(4,1,1)$, $(8,2,1)$, and $(11,1,4)$.

If $(\hat{q}, e, s_3)=(8,2,1)$, then $s_1=1$.
This contradicts Lemma \ref{lemma-torsion-d}.

Let $(\hat{q}, e, s_3)=(4,1,1)$. Then
$s_1=0$ and $\Cl(\hat{X})\simeq\ZZ$ (again by Lemma \ref{lemma-torsion-d}).
In particular, $\bar{f}$ is birational and $\dim|\Theta|\ge3$. By Lemma \ref{lemma-Fano-dimA-ge3} we have
$\hat{X}\simeq\PP^3$ and $|\Theta|=\bar{f}_* \bar\MMM_3$.
Then from \eqref {equation-b-gamma-delta-1} we obtain
\begin{equation} \label{equation-no-q=7-B2334b-gamma-delta-1-l}
\begin{array}{lll}
b &=& 4 \delta-7,
\\[1pt]
\gamma_3 &=& \delta-3.
\end{array}
\end{equation}
Since the linear system $\bar{f}_* \bar\MMM_3=|\Theta|$
has no base points, we have $\gamma_3=0$, $\delta=3$, and $b=5$.
Therefore, $\hat{P}:=\bar{f}(\bar F)$ is a point.
Consider the linear subsystem $\hat\LLL \subset |\Theta|$
consisting of all the divisors passing through $\hat{P}$.
Then
$\bar\MMM_3=\bar{f}^*\hat\LLL=\bar\LLL+a \bar F$,
where $\bar\LLL $ is linear system without
fixed component and $a \in \ZZ_{>0}$.
Hence, $|3A|=\MMM_3=\LLL+a F$.
Then $\LLL \subset|k A|$, where $k \le 2$.
Since $\dim \LLL\ge\dim \hat\LLL=2$, this is impossible.

Finally, let $(\hat{q}, e, s_3)=(11,1,4)$.
Then $s_1=1=e$. Since in this case $\dim|\Theta|\le 0$ (see \cite[Prop. 3.6]{Prokhorov2008a}), 
we have $\bar{f}(\bar E)=\bar{f}(\bar M_1)$ which is an absurd.
The contradiction completes the proof of Proposition \ref{proposition-no-q=7-B2334-c}.
\end{proof}

\begin{mtheorem}{\bf Proposition.} \label{proposition-no-q=7-B2334-dp}
A general member $M \in | 3A|$ is a del Pezzo surface of degree
$4$. Its singular locus consists of two
point of types $\mathrm {A_1}$ and $\mathrm {A_3}$ or a single point of type $\mathrm {D_5}$.
\end{mtheorem}
\begin{proof}
According to Proposition \ref{proposition-no-q=7-B2334-c} the
pair $(X, \MMM)$ has only canonical singularities.
Then a general element of $M \in \MMM$
is a normal surface with only Du Val singularities
\cite[1.12, 1.21]{Alexeev-1994ge}. By the adjunction formula
$-K_M=4A|_M$ and $K_M^2=4$.
Therefore, $M$ is del Pezzo surface of degree $4$.

Let $l:=M \cap M_1$, where $M_1 \in | A|$.
Since $3 M_1 \in \MMM_3$, we have $l \supset \Bs \MMM_3$.
By Bertini's theorem the surface $M$ is smooth outside of $l$.
Since $-K_M \cdot l=4A\cdot M \cdot A=1$, the curve $l$
is reduced irreducible.
Apply the adjunction formula
(see \cite[ch 16.]{Utah}):
\[
K_{l}+\Diff_{l} (0)=(K_{M}+l) |_{l} =-3A|_{L}.
\]
Hence $\deg K_l=-2$, $l\simeq\PP^1$, and $\deg \Diff_{l} (0)=5/4$.

First, consider the case where the point $P_4 \in X$ is a
cyclic quotient. Then $X$ has also a cyclic quotient singularity $P_2$
of type $\frac12 (1,1,1)$.
Since both $M$ and $M_1$ pass through $P_2$ and $P_4$, we have $P_2, \, P_4 \in \Supp (\Diff_{l} (0))$.
We write $\Diff_{l} (0)=\lambda_2 P_2+\lambda_4 P_4+D$,
where $D$ is effective $\QQ$-divisor.
Note that coefficients of the different lying in the range $(0,1)$,
have the form $1-1/m_i$ for some integer $m_i$ (see \cite[ch. 16]{Utah}).
Taking this into account we get a unique possibility: $\Diff_{l} (0)=\frac 12 P_2+\frac 34 P_4$.
This means that
$\Sing(M)=\{P_2, \, P_4\}$ and the singularity of $M$ at $P_2$ (resp. $P_4$)
is of type $\mathrm {A_1}$ (resp. $\mathrm {A_3}$).

Suppose now that the point $P_4 \in X$ of index $4$ is of type $\mathrm {cAx/4}$.
Since $M \sim -K_X$ in a neighborhood of $P_4$, the point $M \ni P_4$ cannot be
Du Val of type $\mathrm {A_n}$ nor $\mathrm {D_4}$ (see \cite[6.4B]{Reid-YPG1987}).
According to Noether's formula applied to the minimum resolution,
the point $M \ni P_4$ is of type $\mathrm {D_5}$ and
$M$ has no other singularities.
\end{proof}

\begin{scase}{\bf Example.} \label{exampleDu-Val-surfaceA1A3D5}
Let $M\subset \PP (1,2,3,4)$ be  the surface 
given by one of the following two equations of degree $6$:
\begin{gather}
\label{Du-Val-surfaceA1A3}
x_2x_4+x_3^2+x_1^6=0,
\\
\label{Du-Val-surfaceD5}
x_1^2x_4+x_3^2+x_2^3=0,
\end{gather}
where $x_k$ are homogeneous coordinates with $\deg x_k=k$.
Then $M$ is del Pezzo surface of degree $4$ with Du Val singularities.
The  singular locus of $M$ is consists of two points of  types $\mathrm {A_1}$ and $\mathrm {A_3}$ in the case 
\eqref{Du-Val-surfaceA1A3} and one point of type $\mathrm {D_5}$ in the case 
\eqref{Du-Val-surfaceD5}. The curve $l:=\{x_1=0\} \cap M$ is irreducible
and it is a line in the anticanonical embedding.
Note that the equality $-K_M \sim 4l$ holds.
\end{scase}

\begin{scase}{\bf Remark.} \label{remarkDu-Val-surfaceA1A3D5}
It is easy to show that any surface of degree $6$ in
$\PP (1,2,3,4)$ with only Du Val singularities up to a coordinate change 
has the form \eqref {Du-Val-surfaceA1A3} or
\eqref{Du-Val-surfaceD5}.
\end{scase}

\begin{mtheorem}{\bf Lemma.} \label{lemma-DuVal-A1A3D5}
Let $M$ be a del Pezzo surface of degree $4$ whose singular locus consists of
\begin{enumerate} [leftmargin=30pt]
\renewcommand\theenumi {{\rm \alph {enumi})}}
\renewcommand\labelenumi {{\rm \alph {enumi})}}
\item
two points of types $\mathrm {A_1}$ and $\mathrm {A_3}$, or
\item
one point of type $\mathrm {D_5}$.
\end{enumerate}
Then $M$ is unique up to isomorphism and
admits an embedding $M \subset \PP (1,2,3,4)$ as a surface of degree $6$
defined by the equation \eqref {Du-Val-surfaceA1A3} or
\eqref{Du-Val-surfaceD5}.
\end{mtheorem}

\begin{proof}
The anticanonical model of $M$ is an intersection of two quadrics in $\PP^4$
which is unique up to projective equivalence
(see, e.g., \cite[\S 8.6.1]{Dolgachev-ClassicalAlgGeom}).
According to Example \xref{exampleDu-Val-surfaceA1A3D5} the surface $M$
admits the desired embedding.
\end{proof}

\begin{proof} [Proof of Theorem \xref{proposition-q=7main}
in the case \xref{no-q=7-B2334}]
According to Lemma \xref{lemma-DuVal-A1A3D5} a general surface
$M \in | 3A|$ can be embedded to $\PP (1,2,3,4)$ and is defined there by a
quasihomogeneous polynomial $s$ of the form \eqref {Du-Val-surfaceA1A3} or
\eqref{Du-Val-surfaceD5}.
Therefore, there is an isomorphism of graded algebras
\[
R(M, l):=\bigoplus_{n\ge0} H^0(M, \OOO_M (nl))\simeq\Bbbk [x_1, x_2, x_3, x_4]/(s).
\]
Further, as in \cite[\S 3]{Prokhorov-2013-fano}, we employ the standard hyperplane section principle
(see, e.g., \cite{Mori-1975}).
Since $H^1 (X, \OOO_X (mA))=0$ for $m>-7$, the natural restriction homomorphism
\[
\Phi: R(X, A):=\bigoplus_{n\ge0} H^0(X, \OOO_X (nA)) \longrightarrow R(M, l)
\]
is surjective.
Let $y_i \in R(X, A)$ be any elements such that $\Phi (y_i)=x_i$.
The kernel of $\Phi$ is a principal ideal generated by a
homogeneous element $y_3'$ of degree $3$ (the equation of $M$).
Thus, $R(M, l)=R(X, A)/(y_3')$.
Then the elements of $y_1, y_2, y_3, y_3', y_4$ generate the algebra $R(X, A)$
and
\[
R(X, A)\simeq\Bbbk [y_1, y_2, y_3, y_3', y_4]/(v),
\]
where
$v \in \Bbbk [y_1, y_2, y_3, y_3', y_4] $ is a homogeneous element such that $v \equiv s \mod y_3'$.
This proves Theorem \xref{proposition-q=7main} in the case \xref{no-q=7-B2334}.
\end{proof}

 
  \newcommand{\etalchar}[1]{$^{#1}$}
\def\cprime{$'$}

\end{document}